# Approximate Series Solution of Nonlinear, Fractional Klein-Gordon Equations Using Fractional Reduced Differential Transform Method

[1]Eman Abuteen, [2]Asad Freihat, [2]Mohammed Al-Smadi, [3]Hammad Khalil and [3]Rahmat Ali Khan

[1]Department of Applied Science, Faculty of Engineering Technology, Al-Balqa Applied University, Amman 11942, Jordan
[2]Department of Applied Science, Ajloun College, Al-Balqa Applied University, Ajloun 26816, Jordan
[3]Department Mathematics, University of Malakand, Chakadara Dir (L), Khyber Pakhtunkhwa, Pakistan

**Abstract:** This analysis proposes an analytical-numerical approach for providing solutions of a class of nonlinear fractional Klein-Gordon equation subjected to appropriate initial conditions in Caputo sense by using the Fractional Reduced Differential Transform Method (FRDTM). This technique provides the solutions very accurately and efficiently in convergent series formula with easily computable coefficients. The behavior of the approximate series solution for different values of fractional-order $\alpha$ is shown graphically. A comparative study is presented between the FRDTM and Implicit Runge-Kutta approach to illustrate the efficiency and reliability of the proposed technique. Our numerical investigations indicate that the FRDTM is simple, powerful mathematical tool and fully compatible with the complexity of such problems.

## Introduction

Fractional Partial Differential Equations (FPDEs) are widely used in interpretation and modeling of many of realism matters appear in applied mathematics and physics including fluid mechanics, electrical circuits, diffusion, damping laws, relaxation processes, mathematical biology (Klimek, 2005; Kilbas *et al*., 2010; Baleanu *et al*., 2009; Jumarie, 2009; Ortigueira, 2010; Mainardi, 2010). Fractional derivatives provide more accurate models of real-world problems than integer-order derivatives; they are actually found to be a powerful tool to describe certain physical problems. The topic of fractional calculus is a significantly important, useful branch of mathematics, plays a critical and serious role to describe a complex dynamical behavior in tremendous scope of application fields, helps to understand the nature of matter as well as simplified the controlling design without any loss of hereditary behaviors and explain even more complex structures.

Consider the following nonlinear Klein-Gordon equations of one-dimensional time fractional model:

$$\frac{\partial^{2\alpha} u(x,t)}{\partial t^{2\alpha}} = \frac{\partial^2 u(x,t)}{\partial x^2} + au(x,t) + bG(u(x,t)) = f(x,t), t \geq 0 \qquad (1)$$

with initial conditions

$$u(x,0) = g_0(x) \frac{\partial}{\partial t} u(x,0) = g_1(x)$$

where, $a$ and $b$ are real constants, $f(x,t)$, $g_0(x)$ and $g_1(x)$ are known analytical functions, $G(u)$ is a nonlinear function, $u$ is an unknown function of $x$ and $t$ to be determined. This model is derived from well-known Klein-Gordon equations (KGEs) by replacing the time order derivative with fractional derivative of order $\alpha$.

The KGEs are fundamental class of nonlinear evolution equations arising in classical relativistic and quantum mechanics. It got a lot of attention for studying solitons and condensed matter physics (Yusufoglu, 2008; Sweilam *et al*., 2012). On the other hand, analytical-numerical studies of the solution for the FKGEs with the Caputo or Riemann-Liouville fractional derivative were used to handle these problems (Golmankhaneh and Baleanu, 2011; Gepreel and Mohamed, 2013; Yang *et al*., 2014; Marasi and Karimi, 2014; Khader *et al*., 2014). As we know, there is no classical method to handle the nonlinear FPDEs and provide its explicit solution due to the complexities of fractional calculus involving these equations. For this reason, we need a reliable numerical approach to find the coefficients of the fractional series solutions of such equations. During the past few decades, many numerical-analytical methods were developed for handling the FPDEs and their system. For examples of these methods, we refer to the work in (Hesameddini and Fotros, 2012; Moaddy *et al*., 2011; Abdulaziz *et al*., 2008; Hashim *et al*., 2009; Odibat and Momani, 2006; Khalil *et al*., 2015a; 2015b; 2015c; El-Ajou *et al*., 2015;

---

Abu-Gdairi *et al.*, 2015; Freihat and Al-Smadi, 2013; Momani *et al.*, 2014; Al-Smadi *et al.*, 2013; 2015; 2016). On the other hand, many applications for different problems by using other numerical algorithms can be found in (Abu Arqub *et al.*, 2012; 2013; 2014; 2015; Moaddy *et al.*, 2015; Komashynska *et al.*, 2014; 2016).

In this analysis, we intend the application of FRDTM to provide numerical analytical solutions for a class of nonlinear partial differential equations included some well-known fractional Klein-Gordon equations. The FRDTM has several advantages for dealing directly with suggested equations; it needs a few iterations to get high accuracy, it is very simple for obtaining analytical-approximate solutions in rapidly convergent formulas, it allows better significantly information in providing continuous representation of these approximations, and it has the ability for solving other problems appearing in several scientific fields.

This article is organized as follows: in the next section, necessary details and preliminaries about the fractional calculus theory are briefly provided. In section 3, the procedure of the RDTM is presented to construct and predict the series solution for fractional PDEs (1). In section 4, numeric results for certain types of FKGEs are given to verify the validity and performance of the present method. Finally, this article ends with some concluding remarks.

## Mathematical Preliminaries

The basic preliminaries, concepts and notations of fractional integrals and derivatives in Caputo definition (Caputo, 1967) are introduced as follows. Here, we adopt the Caputo fractional derivative, which is a modification of Riemann-Liouville, because the initial conditions that defined during the formulation of the system are similar to those conventional conditions of integer order.

### Definition 1

A real function $u(x,t)$, $x \in \mathbb{R}$, $t > 0$ is said to be in the space $C_\mu$, $\mu \in \mathbb{R}$, if there exists a real number $q > \mu$ such that $u(x,t) = t^q u_1(x,t)$, where $u_1(x,t) \in C(\mathbb{R} \times [0,\infty))$ and it is said to be in the space $C_\mu^m$ if $\frac{\partial^m}{\partial t^m} u(x,t) \in C_\mu, m \in \mathbb{N}$.

### Definition 2

The Riemann-Liouville integral operator of order $\alpha \geq 0$ of a function $u(x,t) \in C\mu$, $\mu \geq -1$ is defined as:

$$J_t^\alpha u(x,t) = \begin{cases} \frac{1}{\Gamma(\alpha)} \int_0^t (t-\xi)^{\alpha-1} u(x,\xi) d\xi, & \alpha > 0, 0 < \xi < t \\ u(x,t), & \alpha = 0 \end{cases} \quad (2)$$

Consequently, the operator $J_t^\alpha$ has the following properties: For $u(x,t) \in C_\mu, \mu \geq -1$, $\alpha, \beta \geq 0$, $c \in \mathbb{R}$ and $\gamma > -1$, one can get:

- $J_t^\alpha J_t^\beta u(x,t) = J_t^{\alpha+\beta} u(x,t) = J_t^\beta J_t^\alpha u(x,t)$
- $J_t^\alpha c = \frac{c}{\Gamma(\alpha+1)} t^\alpha$
- $J_t^\alpha t^\gamma = \frac{\Gamma(\gamma+1)}{\Gamma(\alpha+\gamma+1)} t^{\alpha+\gamma}$

Now, we introduce a modified fractional differential operator $D_t^\alpha$ proposed by Caputo as follows:

$$D_t^\alpha f(x) = J_t^{m-\alpha}(x) = \frac{1}{\Gamma(m-\alpha)} \int_t^x (x-\eta)^{m-\alpha-1} f^{(m)}(\eta) d\eta, t \geq 0$$

For $m-1 < \alpha \leq m$, $m \in \mathbb{N}$, $x \geq t$ and $f x \in C_{-1}^m$.

### Definition 3

For $m$ to be the smallest integer that exceeds $\alpha$, the Caputo time-fractional derivative operator of order $\alpha > 0$ is defined as:

$$D_t^\alpha u(x,t) = \frac{\partial^\alpha u(x,t)}{\partial t^\alpha} = \begin{cases} J_t^{m-\alpha}\left(\frac{\partial^m u(x,t)}{\partial t^m}\right), & 0 \leq m-1 < \alpha < m \\ \frac{\partial^m u(x,t)}{\partial t^m}, & \alpha = m \in \mathbb{N} \end{cases} \quad (3)$$

### Theorem 1

If $m-1 < \alpha \leq m$, $m \in \mathbb{N}$, $u(x,t) \in C_\gamma^m$ and $\gamma \geq -1$, then $D_t^\alpha J_t^\alpha u(x,t) = u(x,t)$ and $J_t^\alpha D_t^\alpha u(x,t) = u(x,t) - \sum_{k=0}^{m-1} \frac{\partial^k u(x,0^+)}{\partial t^k} \frac{t^k}{k!}$, where $t > 0$.

For more details about FDEs, see (Millar and Ross, 1993; Podlubny, 1999; Samko *et al.*, 1993).

## Description of the method

Let $u(x,t)$ be a function of two variables such that $u(x,t) = f(x)g(t)$, then from the properties of the one-dimensional Differential Transform (DT) method, we have:

$$u(x,t) = \sum_{i=0}^{\infty} f(i)x^i \sum_{j=0}^{\infty} g(j)t^j = \sum_{i=0}^{\infty}\sum_{j=0}^{\infty} U(i,j)x^i t^j \quad (4)$$

where, $U(i,j) = f(i)g(j)$ is called the spectrum of $u(x,t)$.

Next, we assume that $u(x,t)$ is continuously differentiable function with respect to space variable $x$ and time $t$.

### Lemma 1 (Srivastava et al., 2013)

Let $u(x,t)$ be an analytic function, then the FRDT of $u$ is given by:

$$U_k(x) = \frac{1}{\Gamma(k\alpha+1)}\left[\frac{\partial^{k\alpha} u(x,t)}{\partial t^{k\alpha}}\right]_{t=t_0} \quad (5)$$

where, $\alpha$ is a parameter describing the order of time-fractional derivative in Caputo sense.

The inverse transformed of $U_k$ is defined by:

$$u(x,t) = \sum_{k=0}^{\infty} U_k(x)(t-t_0)^{k\alpha} \quad (6)$$

From Equations 5 and 6, we have that

$$u(x,t) = \sum_{k=0}^{\infty}\frac{1}{\Gamma(k\alpha+1)}\left[\frac{\partial^{k\alpha} u(x,t)}{\partial t^{k\alpha}}\right]_{t=t_0}(t-t_0)^{k\alpha} \quad (7)$$

In particular, for $t = 0$, Equation 7 reduces to:

$$u(x,t) = \sum_{k=0}^{\infty}\frac{1}{\Gamma(k\alpha+1)}\left[\frac{\partial^{k\alpha} u(x,t)}{\partial t^{k\alpha}}\right]_{t=t_0} t^{k\alpha} \quad (8)$$

Moreover, if $\alpha = 1$, then the FRDT of Equation 7 reduces to the classical RDT method.

From the above lemma, the fundamental operations of the FRDTM are given by the following theorems (Srivastava et al., 2013):

### Theorem 2

Let $u(x,t)$, $v(x,t)$ and $w(x,t)$ be any analytic functions such that $u(x,t) = R_D^{-1}[U_k(x)]$, $v(x,t) = R_D^{-1}[V_k(x)]$ and $w(x,t) = R_D^{-1}[W_k(x)]$, then:

- If $u(x,t) = v(x,t) \pm w(x,t)$, then $U_k(x) = V_k(x) \pm W_k(x)$
- If $u(x,t) = av(x,t)$, then $U_k(x) = aV_k(x)$, $a$ is an arbitrary constant
- If $u(x,t) = x^m t^n v(x,t)$, then $U_k(x) = V_{k-n}(x)$
- If $u(x,t) = x^m t^n$, then $U_k(x) = x^m \delta(k-n)$, $\delta(k) = \begin{cases} 1, k = 0 \\ 0, k \neq 0 \end{cases}$
- If $u(x,t) = v(x,t)\, w(x,t)$, then $U_k(x) = \sum_{r=0}^{k} V_r(x)W_{k-r}(x) = \sum_{r=0}^{k} W_r(x)V_{k-r}(x)$

### Theorem 3

Let $u(x,t)$ and $v(x,t)$ be any two analytic functions such that $u(x,t) = R_D^{-1}[U_k(x)]$ and $v(x,t) = R_D^{-1}[V_k(x)]$, then:

- If $u(x,t) = \dfrac{\partial^r}{\partial x^r} v(x,t)$, then $U_k(x) = \dfrac{\partial^r}{\partial x^r} V_k(x)$
- If $u(x,t) = \dfrac{\partial^{r\partial}}{\partial t^{r\partial}} v(x,t)$, then $U_k(x) = \dfrac{\Gamma(\alpha k + \alpha r + 1)}{\Gamma(\alpha k + 1)} V_{k+r}(x)$

### Corollary 1

If $u(x,t) = e^{\lambda t + \mu x}$, then $U_k(x) = \dfrac{\lambda^k}{k!} e^{\mu x}$.

### Corollary 2

If $u(x,t) = \sin(\eta x + \omega t)$, $v(x,t) = \cos(\eta x + \omega t)$, then $U_k(x) = \dfrac{\omega^k}{k!}\sin\left(\eta x + \dfrac{\pi k}{2}\right)$ and $V_k(x) = \dfrac{\omega^k}{k!}\cos\left(\eta x + \dfrac{\pi k}{2}\right)$.

The reader is referred to (Keskin and Oturanc, 2009; Abazari and Abazari, 2012; Secer, 2012; Sohail and Mohyud-Din, 2012; Al-Amr, 2014) and the references therein to know more details about the reduced differential transform technique, including their applications in various kinds of differential equations.

Now, by applying the FRDTM to Equation 1, we obtain the following recurrence relation formula:

$$\frac{\Gamma(\alpha k + 2\alpha + 1)}{\Gamma(\alpha k + 1)} U_{k+2}(x) = \frac{\partial^2}{\partial x^2} U_k(x) - aU_k(x) - bG[U_k(x)] + F_k(x) \quad (9)$$

where, $F_k(x)$ and $G[U_k(x)]$ are the reduced transformation of the functions $f(x,t)$ and $G(u(x,t))$, respectively.

Using the initial conditions, we have:

$$U_0(x) = g_0(x), U_1(x) = g_1(x) \quad (10)$$

Substituting Equation 10 into Equation 9 and by straightforward iterative calculation, we obtain the values $U_k(x)$, for $k = 1,2,3,\ldots$. Thus, inverse RDT of $\{U_k(x)\}_{k=1}^n$ yields that

$$\breve{u}_n(x,t) = \sum_{k=0}^n U_k(x)(t-t_0)^{k\alpha} \tag{11}$$

Therefore, the closed form solution given by:

$$u(x,t) = \lim_{x \to \infty} \breve{u}_n(x,t) \tag{12}$$

## Numerical Examples

In this section, some numerical examples are given to verify the simplicity and applicability of the present technique in finding approximate series solution for fractional KGEs. The simulation results indicate that the FRDT method is highly accurate and fully compatible with the complexity of the PDEs of fractional-order.

*Example 4.1*

We consider the following one-dimensional linear fractional Klein-Gordon equation:

$$\frac{\partial^\alpha u(x,t)}{\partial t^\alpha} - \frac{\partial^2 u(x,t)}{\partial x^2} - u(x,t) = 0, 0 < \alpha \leq 1 \tag{13}$$

with initial condition:

$$u(x,0) = 1 + \sin(x)$$

where, $x,t \geq 0$.

Applying the transformation (Li and He, 2010), we have:

$$\frac{\partial u(x,t)}{\partial T} - \frac{\partial^2 u(x,t)}{\partial x^2} = u(x,t)$$

By using the FRDTM of Equation 13, we have the recurrence relation formula:

$$U_{k+1}(x) = \frac{\Gamma(\alpha k+1)}{\Gamma(\alpha(k+1)+1)}\left[\frac{\partial^2}{\partial x^2}U_k(x) + U_k(x)\right] \tag{14}$$

with transformed initial data:

$$U_0(x) = 1 + \sin(x) \tag{15}$$

Substituting the condition (15) into Equation 14, we get the values successively $U_k(x)$, $k = 1,2,3,\ldots$, as follows:

$$U_1(x) = \frac{1}{\Gamma(\alpha+1)}, U_2(x) = \frac{1}{\Gamma(2\alpha+1)},$$

$$U_3(x) = \frac{1}{\Gamma(3\alpha+1)},\ldots,U_k(x) = \frac{1}{\Gamma(k\alpha+1)}$$

Thus, the approximate solution can be obtained by:

$$u(x,t) = 1 + \sin(x) + \frac{1}{\Gamma(\alpha+1)}T$$
$$+ \frac{1}{\Gamma(2\alpha+1)}T^2 + \frac{1}{\Gamma(3\alpha+1)}T^3 + \ldots$$

The inverse RDTM is given by:

$$u(x,t) = 1 + \sin(x) + \frac{1}{\Gamma(\alpha+1)}t^\alpha$$
$$+ \frac{1}{\Gamma(2\alpha+1)}t^{2\alpha} + \frac{1}{\Gamma(3\alpha+1)}t^{3\alpha} + \ldots$$
$$= \sin(x) + \sum_{k=0}^\infty \frac{1}{\Gamma(2\alpha+1)}t^{k\alpha}$$

Consequently, the reduced inverse transformed of $U_k(x)$ follows the closed form solution. Setting $\alpha = 1$, the exact solution is:

$$u(x,t) = \sin(x) + e^t$$

To demonstrate the efficiency of the present method, we compare the FRDT approximation with the Implicit Runge-Kutta (IRK) method for $\alpha = 1$. Figure 1 shows the phase portrait of solutions for Example 4.1 using the FRDTM and IRKM for $t \in [0,0.8]$ and $x \in [0,4]$. The numerical results for different time levels of $\alpha$ are presented in Fig. 2. Here, we note that the approximate FRDT are efficiency for time-fractional KGE and very closed to the IRK solutions.

*Example 4.2*

We consider the nonlinear Klein-Gordon fractional model in the form:

$$\frac{\partial^\alpha u(x,t)}{\partial t^\alpha} - \frac{\partial^2 u(x,t)}{\partial x^2} + u^2(x,t) = 0, 0 < \alpha \leq 1 \tag{16}$$

with initial condition:

$$u(x,0) = 1 + \sin(x)$$

Applying the transformation (Li and He, 2010), we have:

$$\frac{\partial u(x,t)}{\partial T} - \frac{\partial^2 u(x,t)}{\partial x^2} = u^2(x,t)$$

By using the FRDTM, we have the recurrence relation formula:

$$U_{k+1}(x) = \frac{\Gamma(\alpha k+1)}{\Gamma(\alpha(k+1)+1)}\left[\frac{\partial^2}{\partial x^2}U_k(x) + \sum_{r=0}^{k}U_r(x)U_{k-r}(x)\right] \quad (17)$$

With RDTM of initial condition:

$$U_0(x) = 1 + \sin(x)$$

Following recurrence relation (17), the sequences components $U_k(x)$, $k = 1,2,3,\ldots$, were computed using the Mathematica package, can be successively given by:

$$U_1(x) = \frac{-\Gamma(\alpha)}{\Gamma(2\alpha)}\left[1 + 3\sin(x) + \sin^2(x)\right],$$

$$U_2(x) = \frac{-\Gamma(\alpha)}{2\Gamma(3\alpha)}\left[\begin{array}{l}12\cos(2x) - 25\sin(x) \\ +\sin(3x) - 12\end{array}\right],$$

$$U_3(x) = \frac{-\Gamma(\alpha)}{8(\Gamma(2\alpha))^2\Gamma(4\alpha)}\left[\begin{array}{l}2\Gamma(\alpha)\Gamma(3\alpha) \\ \left((-3+\cos(2x)-6\sin^2(x))+4(\Gamma(2\alpha))^2\right) \\ \left(\begin{array}{l}49-98\cos(2x)\\ +\cos(4x)+111\sin(x)-23\sin(3x)\end{array}\right)\end{array}\right],$$

⋮

*and so on*

The approximate form solution is given by:

$$u(x,t) = 1 + \sin(x) + U_1(x)T + U_2(x)T^2 + U_3(x)T^3 + \ldots$$

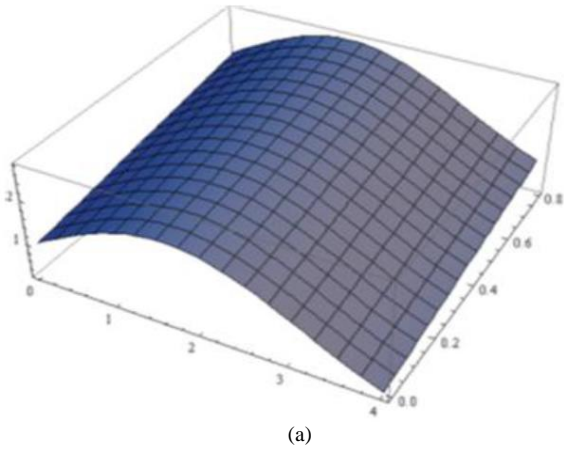
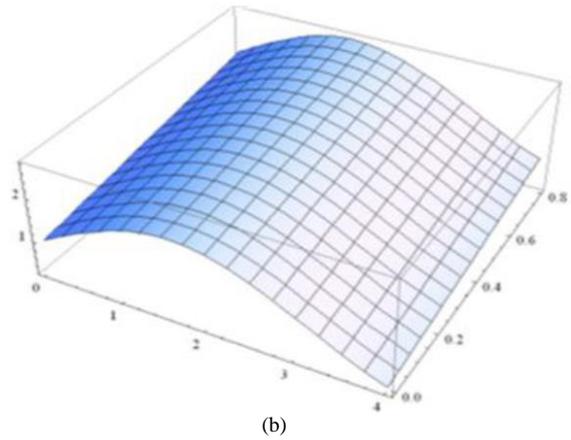

(a)          (b)

Fig. 1. Comparison of phase plot for $u(x,t)$ of Example 4.1 at $\alpha = 1$, $x \in [0,4]$ and $t \in 0,0.8$: (a) the FRDTM; (b) implicit Runge-Kutta method

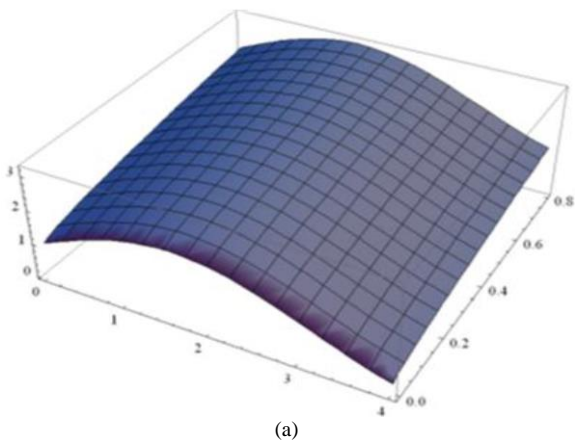
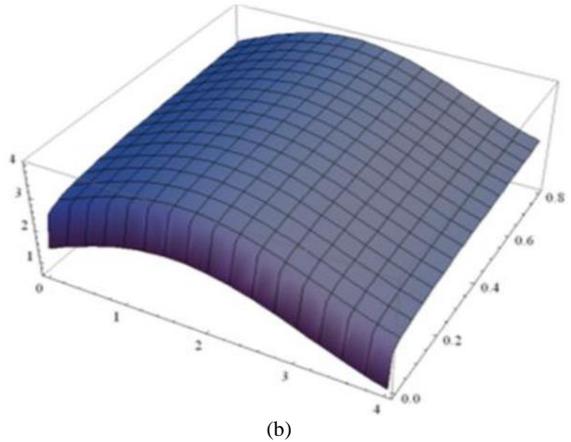

(a)          (b)

Fig. 2. Phase plot of the solution $u(x,t)$ of Example 4.1 using the FRDTM for $x \in [0,4]$ and $t \in [0,0.8]$: (a) at $\alpha = 0.7$; (b) at $\alpha = 0.1$

Accordingly, the inverse RDT is given by:

$$u(x,t) = \sum_{k=0}^{\infty} U_k(x) t^{k\alpha}$$

$$= 1 + \sin(x) - \frac{\Gamma(\alpha)}{\Gamma(2\alpha)}\left[1 + 3\sin(x) + \sin^2(x)\right]t^\alpha$$

$$- \frac{\Gamma(\alpha)}{2\Gamma(3\alpha)}\left[12\cos(2x) - 25\sin(x) + \sin(3x) - 12\right]t^{2\alpha}$$

$$- \frac{\Gamma(\alpha)}{8(\Gamma(2\alpha))^2 \Gamma(4\alpha)} \begin{bmatrix} 2\Gamma(\alpha)\Gamma(3\alpha)\begin{pmatrix} -3+\cos(2x) \\ -6\sin^2(x) \end{pmatrix} \\ +4(\Gamma(2\alpha))^2 \begin{pmatrix} 49 - 98\cos(2x) + \cos(4x) \\ +111\sin(x) - 23\sin(3x) \end{pmatrix} \end{bmatrix} t^{3\alpha}$$

$$+ \cdots.$$

The series solution of $u$ follows closed form solution.

The geometric behaviors of the solution for Example 4.2 are studied by drawing the 3-dimensional space figures of the FRDT approximate solution together with its corresponding IRK solution. Figure 3 shows the comparison between the FRDT approximate solution and IRK solution at $\alpha = 1$ for $x \in [-2,2]$ and $t \in [0,0.01]$. Whileas, Fig. 4 shows the solution behavior of the nonlinear FKGE for different specific cases of $\alpha$ in the domain $x \in [-2,2]$ and $t \in [0,0.01]$. The performance errors analysis are obtained by the FRDTM at $x = 2$ and summarized in Table 1. Numerically, it is showed that the RDT method is effective and accurate.

### Example 4.3

We consider the nonlinear Klein-Gordon fractional model in the form:

$$\frac{\partial^\alpha u(x,t)}{\partial t^\alpha} - \frac{\partial^2 u(x,t)}{\partial x^2} + u(x,t) - u^3(x,t) = 0, \quad 0 < \alpha \leq 1 \quad (18)$$

with initial condition:

$$u(x,0) = -\sec h(x)$$

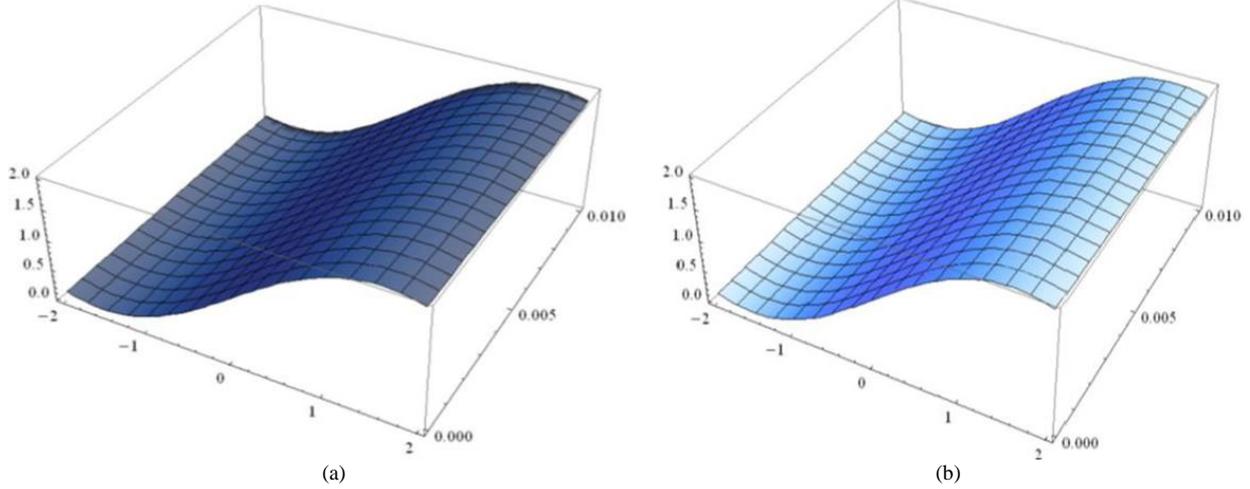

(a)          (b)

Fig. 3. Comparison of phase plot for $u(x,t)$ of Example 4.2 at $\alpha = 1$, $x \in [-2,2]$ and $t \in 0,0.01$ : (a) the FRDTM; (b) implicit RKM

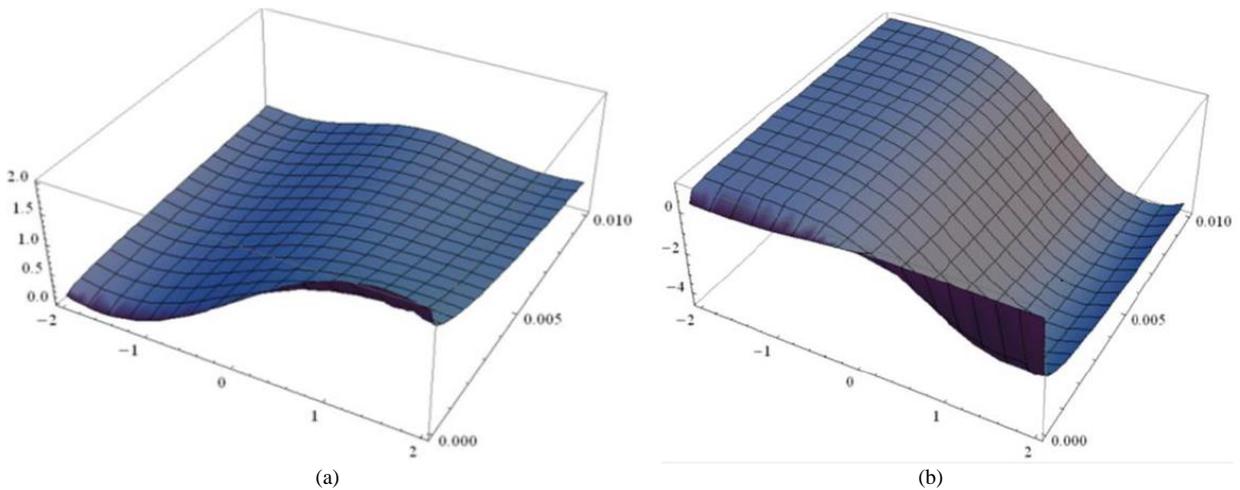

(a)          (b)

Fig. 4. Phase plot of the solution $u(x,t)$ of Example 4.2 using the FRDTM for $x \in [-2,2]$ and $t \in [0,0.01]$: (a) at $\alpha = 0.4$; (b) at $\alpha = 0.1$

Table 1. The error analysis for Example 4.2 when $\alpha = 1$ with $x = 2$.

| $t$ | RDTM | IRKM | Absolute error | Relative error |
|---|---|---|---|---|
| 0.000 | 1.909297426825682 | 1.909297426825682 | 0.00000 | 0.00000 |
| 0.001 | 1.904742712734773 | 1.904743294464910 | $5.81730 \times 10^{-7}$ | $3.05411 \times 10^{-7}$ |
| 0.002 | 1.900187998643864 | 1.900189162104138 | $1.16346 \times 10^{-6}$ | $6.12287 \times 10^{-7}$ |
| 0.003 | 1.895633284552955 | 1.895635029743366 | $1.74519 \times 10^{-6}$ | $9.20636 \times 10^{-7}$ |
| 0.004 | 1.891078570462046 | 1.891080897382594 | $2.32692 \times 10^{-6}$ | $1.23047 \times 10^{-6}$ |
| 0.005 | 1.886523856371138 | 1.886526765021823 | $2.90865 \times 10^{-6}$ | $1.54180 \times 10^{-6}$ |
| 0.006 | 1.881969142280229 | 1.881972632661051 | $3.49038 \times 10^{-6}$ | $1.85464 \times 10^{-6}$ |
| 0.007 | 1.877414428189319 | 1.877418500300279 | $4.07211 \times 10^{-6}$ | $2.16899 \times 10^{-6}$ |
| 0.008 | 1.872859714098411 | 1.872864367939507 | $4.65384 \times 10^{-6}$ | $2.48488 \times 10^{-6}$ |
| 0.009 | 1.868305000007502 | 1.868310235578735 | $5.23557 \times 10^{-6}$ | $2.80230 \times 10^{-6}$ |
| 0.010 | 1.863750285916593 | 1.863756103217964 | $5.81730 \times 10^{-6}$ | $3.12128 \times 10^{-6}$ |

Applying the transformation (Li and He, 2010), we have:

$$\frac{\partial u(x,t)}{\partial T} = \frac{\partial^2 u(x,t)}{\partial x^2} - u(x,t) + u^3(x,t)$$

By using the FRDTM, we have the recurrence relation formula:

$$U_{k+1}(x) = \frac{\Gamma(\alpha k + 1)}{\Gamma(\alpha(k+1)+1)} \left[ \frac{\partial^2}{\partial x^2} U_k(x) - U_k(x) + \sum_{r=0}^{k} \sum_{i=0}^{r} U_i(x) U_{k-r}(x) \right] \quad (19)$$

With RDT of initial condition:

$$U_0(x) = -\sec h(x)$$

Following recurrence relation (19), the sequences components $U_k(x)$, $k = 1,2,3,\ldots$, were computed using the Mathematica 9.0 package, can be successively given by:

$$U_1(x) = \frac{\Gamma(\alpha)}{\Gamma(2\alpha)} \left[ (\sec h(x))^3 \right],$$

$$U_2(x) = \frac{\Gamma(\alpha)}{\Gamma(3\alpha)} \left[ (\sec h(x))^5 (4\cosh(2x) - 5) \right],$$

$$U_3(x) = \frac{\Gamma(\alpha)}{(\Gamma(2\alpha))^2 \Gamma(4\alpha)} \left[ (\sec h(x))^7 \begin{pmatrix} (123 - 112\cosh(2x)) \\ +8\cosh(4x) \\ (\Gamma(2\alpha))^2 3\Gamma(\alpha)\Gamma(3\alpha) \end{pmatrix} \right],$$

⋮

*and so on*

By taking the inverse RDT of $\{U_k(x)\}_{k=0}^{n}$, we have

$$u(x,t) = \sum_{k=0}^{\infty} U_k(x) t^{k\alpha}$$

$$= -\sec h(x) + \frac{\Gamma(\alpha)}{\Gamma(2\alpha)} (\sec h(x))^3 t^{\alpha}$$

$$+ \frac{\Gamma(\alpha)}{\Gamma(3\alpha)} \left[ (\sec h(x))^5 (4\cosh(2x) - 5) \right] t^{2\alpha}$$

$$+ \frac{\Gamma(\alpha)}{(\Gamma(2\alpha))^2 \Gamma(4\alpha)} \left[ (\sec h(x))^7 \begin{pmatrix} (123 - 112\cosh(2x)) \\ +8\cosh(4x) \\ (\Gamma(2\alpha))^2 - 3(\alpha)\Gamma(3\alpha) \end{pmatrix} \right] t^{3\alpha} + \cdots.$$

Figure 5 and 6 show the solution behavior for the nonlinear time-fractional KGE in Example 4.3 by using the RDTM and IRKM at $x \in [-2,2]$ and $t \in 0,0.01$ for different particular values of $\alpha$.

*Example 4.4*

We consider the nonlinear Klein-Gordon fractional model in the from:

$$\frac{\partial^{2\alpha} u(x,t)}{\partial t^{2\alpha}} - \frac{\partial^2 u(x,t)}{\partial x^2} + u^2(x,t) = 0, 1 < \alpha \leq 2 \quad (20)$$

with initial conditions:

$$u(x,0) = 1 + \sin(x), \frac{\partial}{\partial t} u(x,0) = 0 \quad (21)$$

Using the transformation (Li and He, 2010), we have:

$$\frac{\partial u(x,t)}{\partial T} = \frac{\partial^2 u(x,t)}{\partial x^2} - u^2(x,t)$$

By using the FRDTM, we obtain the recurrence relations formula as:

$$U_{k+2}(x) = \frac{\Gamma(\alpha k + 1)}{\Gamma(\alpha(k+2)+1)} \left[ \frac{\partial^2}{\partial x^2} U_k(x) - \sum_{r=0}^{k} U_r(x) U_{k-r}(x) \right] \quad (22)$$

With RDT of initial condition:

$$U_0(x) = 1 + \sin(x), U_1(x) = 0$$

Following recurrence relation (22) and by straightforward iterative steps, yields:

$$U_2(x) = -\frac{\Gamma(\alpha)}{\Gamma(3\alpha)}\left[1 + 3\sin(x) + \sin^2(x)\right], U_3(x) = 0,$$

$$U_4(x) = -\frac{\Gamma(\alpha)}{2\Gamma(4\alpha)}\left[12\cos(2x) - 25\sin(x) + \sin(3x) - 12\right], U_5(x) = 0,$$

$$U_6(x) = -\frac{\Gamma(\alpha)}{8(\Gamma(2\alpha))^2 \Gamma(5\alpha)} \left[ 2\Gamma(\alpha)\Gamma(3\alpha)\begin{pmatrix}-3+\cos 2x \\ -6\sin^2 x\end{pmatrix} + 4(\Gamma(2\alpha))^2 \begin{pmatrix}-12+12\cos 2x \\ -25\sin x + \sin 3x\end{pmatrix} \right],$$

$$U_7(x) = 0$$

$$\vdots$$

*and so on*

The approximate form solution is obtained by:

$$u(x,t) = 1 + \sin(x) + U_1(x)T + U_2(x)T^2 + U_3(x)T^3 + U_4(x)T^4 + \cdots$$

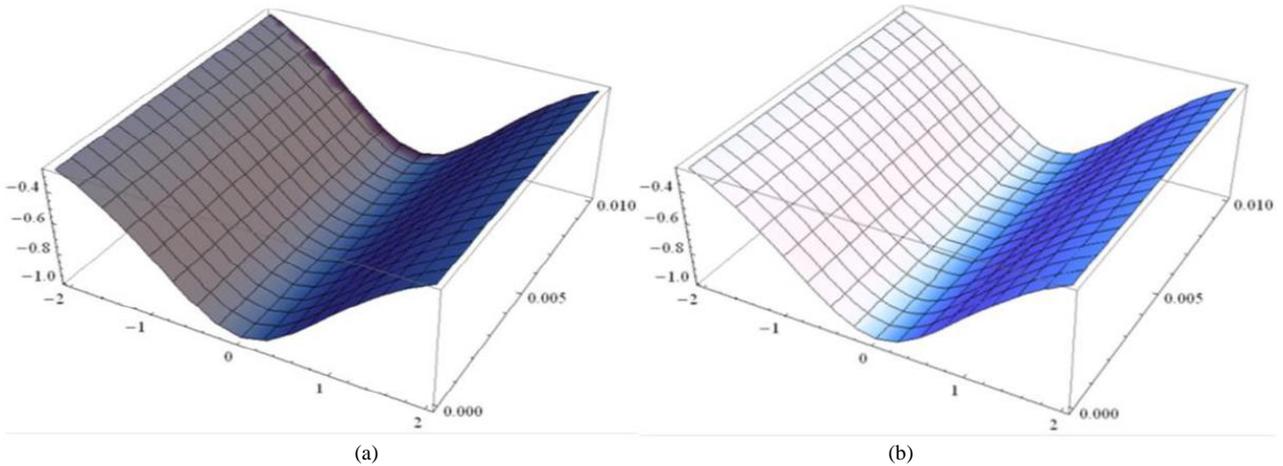

(a)            (b)

Fig. 5. Comparison of phase plot for $u(x,t)$ of Example 4.3 at $\alpha = 1$, $x \in [-2,2]$ and $t \in [0,0.01]$: (a) the FRDTM; (b) implicit Runge-Kutta method

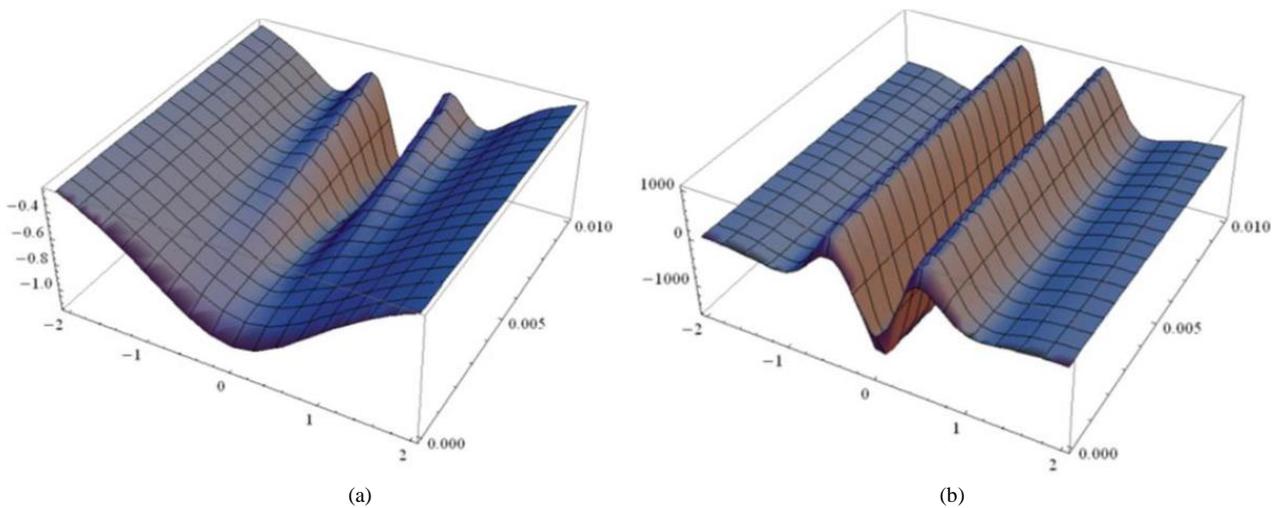

(a)            (b)

Fig. 6. Phase plot of the solution $u(x,t)$ of Example 4.2 using the FRDTM for $x \in [-2,2]$ and $t \in [0,0.01]$: (a) at $\alpha = 0.4$; (b) at $\alpha = 0.1$

Table 2. Numerical results of $u$ when $x = 2$ for Example 4.4

| $t$ | $\alpha = 1.25$ | $\alpha = 1.5$ | $\alpha = 1.75$ | $\alpha = 2$ |
|---|---|---|---|---|
| 0.00 | 1.9092974268257 | 1.9092974268257 | 1.9092974268257 | 1.9092974268257 |
| 0.02 | 1.8577557126584 | 1.8953639720340 | 1.9056740630824 | 1.9083864840075 |
| 0.04 | 1.7867098804114 | 1.8698876653513 | 1.8971099324942 | 1.9056536555530 |
| 0.06 | 1.7057992950300 | 1.8368970719530 | 1.8845189771575 | 1.9010989414620 |
| 0.08 | 1.6177334620125 | 1.7978297884919 | 1.8683037456485 | 1.8947223417348 |
| 0.10 | 1.5239332295153 | 1.7535166664471 | 1.8487201924577 | 1.8865238563711 |
| 0.12 | 1.4252945743833 | 1.7045182992625 | 1.8259529888198 | 1.8765034853711 |
| 0.14 | 1.3224416173263 | 1.6512462328457 | 1.8001449727689 | 1.8646612287348 |
| 0.16 | 1.2158375439572 | 1.5940193350310 | 1.7714116686258 | 1.8509970864620 |
| 0.18 | 1.1058416180037 | 1.5330941474492 | 1.7398494741889 | 1.8355110585530 |
| 0.20 | 0.9927417362082 | 1.4686828986574 | 1.7055407099216 | 1.8182031450075 |

Table 3. The Numerical results for Example 4.4 when $\alpha = 2$ with $x = 2$.

| $t$ | FRDTM | IRKM | Absolute error | Relative error |
|---|---|---|---|---|
| 0.00 | 1.90929742682568 | 1.90929742682568 | 0.00000 | 0.00000 |
| 0.02 | 1.90838648400750 | 1.90838660035353 | $1.16346 \times 10^{-7}$ | $6.09656 \times 10^{-8}$ |
| 0.04 | 1.90565365555296 | 1.90565412093706 | $4.65384 \times 10^{-6}$ | $2.44212 \times 10^{-7}$ |
| 0.06 | 1.90109894146205 | 1.90109998857629 | $1.04711 \times 10^{-6}$ | $5.50794 \times 10^{-7}$ |
| 0.08 | 1.89472234173477 | 1.89472420327121 | $1.86154 \times 10^{-6}$ | $9.82484 \times 10^{-7}$ |
| 0.10 | 1.88652385637114 | 1.88652676502182 | $2.90865 \times 10^{-6}$ | $1.54180 \times 10^{-6}$ |
| 0.12 | 1.87650348537114 | 1.87650767382812 | $4.18846 \times 10^{-6}$ | $2.23205 \times 10^{-6}$ |
| 0.14 | 1.86466122873478 | 1.86466692969012 | $5.70096 \times 10^{-6}$ | $3.05736 \times 10^{-6}$ |
| 0.16 | 1.85099708646205 | 1.85100453260780 | $7.44615 \times 10^{-6}$ | $4.02276 \times 10^{-6}$ |
| 0.18 | 1.83551105855296 | 1.83552048258118 | $9.42403 \times 10^{-6}$ | $5.13425 \times 10^{-6}$ |
| 0.20 | 1.81820314500750 | 1.81821477961024 | $1.16346 \times 10^{-5}$ | $6.39892 \times 10^{-6}$ |

Consequently, the reduced inverse DT of $U_k(x)$ will be given by:

$$u(x,t) = \sum_{k=0}^{\infty} U_k(x) t^{k\alpha}$$
$$= 1 + \sin(x) - \frac{\Gamma(\alpha)}{\Gamma(2\alpha)}[1 + 3\sin(x) + \sin^2(x)] t^{\alpha}$$
$$- \frac{\Gamma(\alpha)}{\Gamma(3\alpha)}[1 + 3\sin(x) + \sin^2(x)] t^{2\alpha}$$
$$- \frac{\Gamma(\alpha)}{2\Gamma(4\alpha)}[12\cos(2x) - 25\sin(x) + \sin(3x) - 12] t^{3\alpha}$$
$$- \frac{\Gamma(\alpha)}{8(\Gamma(2\alpha))^2 \Gamma(5\alpha)} \left[ \begin{array}{c} 2\Gamma(\alpha)\Gamma(3\alpha) \begin{pmatrix} -3 + \cos 2x \\ -6\sin^2 x \end{pmatrix} \\ + 4(\Gamma(2\alpha))^2 \begin{pmatrix} -12 + 12\cos 2x \\ -25\sin x + \sin 3x \end{pmatrix} \end{array} \right] t^{4\alpha} + \cdots$$

In Table 2, we summarized numerical values for the solution of Equation 20 when $\alpha = 1.25, 1.5, 1.75$ and $\alpha = 2$ with $x = 2$. Also numerical results using our approximation and the implicit Runge-Kutta method for $\alpha = 2$ are summarized in Table 3. From these tables, it is easier to observe that the numerical approximations are in agreement with each other's and with the IRK method. Also, they have same behavior as those obtained using the implicit Runge-Kutta method.

# Concluding remarks

In this study, we present numerical algorithm for finding approximate form solutions of a class of Klein-Gordon fractional model based upon FRDTM. This method was used directly without employing linearization and perturbation. The efficiency and capability of the present algorithm have been checked via several illustrated examples. The results reveal the complete reliability of this method with a great potential in scientific applications. Finally, we conclude that the FRDTM is very powerful, straightforward and effective to obtain analytical numerical solutions of a wide variety problems related to fractional PDEs applied in mathematics, physics and engineering. Computations of this paper have been carried out by using the computer package of Mathematica 9.